\documentclass[12pt]{llncs}
\usepackage{amsfonts}

\newcommand{\set}[1]{\left\{ #1\right\}}
\newcommand{\gilt}{:}
\newcommand{\sodass}{\,:\,}
\newcommand{\setGilt}[2]{\left\{ #1\sodass #2\right\}}

\newcommand{\realrange}[2]{\left[#1, #2\right]}

\newcommand{\unitrange}[2]{\realrange{0}{1}}

\newcommand{\llabel}[1]{\label{\labelprefix:#1}}
\newcommand{\labelprefix}{} % later redefined using renewcommand

\newcommand{\discussionsize}{\small}

\newcommand{\frage}[1]{}

\newenvironment{code}{\noindent%\sf%
\begin{tabbing}%
\hspace{2em}\=\hspace{2em}\=\hspace{2em}\=\hspace{2em}\=\hspace{2em}\=%
\hspace{2em}\=\hspace{2em}\=\hspace{2em}\=\hspace{2em}\=\hspace{2em}\=%
\kill}{\end{tabbing}}

\newcommand{\labelcommand}{}
\newcommand{\captiontext}{}
\newsavebox{\codeparam}
\newcounter{lineNumber}
\newenvironment{disscodepos}[3]{%
\renewcommand{\labelcommand}{#2}%
\renewcommand{\captiontext}{#3}%
\sbox{\codeparam}{\parbox{\textwidth}{#3}}%
\begin{figure}[#1]\begin{center}\begin{code}\setcounter{lineNumber}{1}}{%
\end{code}\end{center}\caption{\llabel{\labelcommand}\captiontext}\end{figure}}

{\end{disscodepos}}

\newcommand{\Is}       {:=}

\newdimen\endofsize\endofsize=0.5em
\def\endofbeweis{~\quad\hglue\hsize minus\hsize
                 \hbox{\vrule height \endofsize width
\endofsize}\par}
\usepackage[utf8]{inputenc}
\usepackage{algorithmic}
\usepackage{algorithm}
\usepackage{amsmath}
\usepackage{amssymb}
\usepackage{color}
\usepackage{fullpage}     
\usepackage{latexsym}
\usepackage{makeidx}
\usepackage{multicol}
\usepackage{numprint}
\usepackage{t1enc}
\usepackage{times}
\usepackage{graphicx}
\usepackage{url}
\usepackage{todonotes}
\npdecimalsign{.} % we want . not , in numbers

\definecolor{mygrey}{gray}{0.75}
\newcommand{\ie}{i.e.\ }
\newcommand{\etal}{et~al.\ }
\newcommand{\eg}{e.g.\ }

\def\MdR{\ensuremath{\mathbb{R}}}

\newcommand{\csch}[1]{}
\newcommand{\sout}[1]{}

\newcommand{\mytitle}{Incorporating Road Networks into Territory Design}
\begin{document}
\title{\mytitle}
\institute{}
\author{Nitin Ahuja, Matthias Bender, Peter Sanders, Christian Schulz and Andreas Wagner $^{a,b,c,c,a}$ \\\ \\
$^a$ PTV AG, \{nitin.ahuja, andreas.wagner\}@ptvgroup.com \\
$^b$ FZI Research Centre for Information Technology,  {matthias.bender@fzi.de} \\
$^c$ Karlsruhe Institute of Technology, \{sanders, christian.schulz\}@kit.edu}
\pagestyle{plain}
\maketitle
\begin{abstract}
Given a set of \emph{basic areas}, the territory design problem asks to create a predefined number of territories, each containing at least one basic area, such that an objective function is optimized.
Desired properties of territories often include a reasonable balance, compact form, contiguity and small average journey times which are usually encoded in the objective function or formulated as constraints. 
We address the territory design problem by developing graph theoretic models that also consider the underlying road network.
The derived graph models enable us to tackle the territory design problem by modifying graph partitioning algorithms and mixed integer programming formulations so that the objective of the planning problem is taken into account.
We test and compare the algorithms on several real world instances. 

\end{abstract}
\thispagestyle{empty}
\section{Practical Motivation}

The starting point for this research was a real territory design problem at the PTV AG; design and plan service territories for companies that send agents to customers on a regular basis. In this case, the number of territories to be planned is predetermined by the number of service agents.
Required properties of territories are manifold. For example, the obtained territories should be reasonably balanced with respect to an activity index of the customers, they should have ``good'' accessibility (non-accessible areas have to be accounted for), have to be contiguous and compact; and  most importantly the expected journey time for visiting a typical set of customers within the territory should be small. Even small improvements in the objective save fuel and time that is needed to visit the customers. In turn, the costs of the service company are reduced and the usage of valuable resources is optimized. This optimization can also lead to the company serving even more customers using the same amount of employees or increasing the face time with its customers.

More formally, given a set of \emph{basic areas} such as zip code areas or customer sites, the territory design problem asks to create a predefined number of larger territories such that an objective function encoding the desired properties is optimized and territory constraints are fulfilled. Apart from the introduced example, other important applications include political districting, electrical power restricting or planning territories for social facilities like hospitals or administrative units.

We organize the paper as follows.
We begin in Section~\ref{s:preliminaries} by introducing basic concepts and by summarizing related work. 
The main parts of the paper are Section~\ref{s:territorydesignbygp} and Section~\ref{s:territorydesignbymip}. 
The former introduces our rationale, graph models and the approach employing graph partitioning algorithms to perform territory design. 
The latter uses a mixed integer programming approach to tackle the problem.
We report on carefully designed experiments performed on several real world instances in Section~\ref{s:experiments}.
Experiments indicate that our algorithms compute high quality territories that have small average journey times in particular. 
Finally, we conclude with Section~\ref{s:conclusion}.
\csch{TODO: add MIP}

\section{Preliminaries}
\label{s:preliminaries}
\subsection{Basic concepts}
Consider a set of \emph{basic areas} $B=\{b_1, \ldots, b_n\}$ with an \emph{activity index} $a: B \to \MdR_{\geq 0}$. 
We use the short notation $a_i$ for $a(b_i)$ and extend the activity index to sets, \ie , for a set $B' \subseteq B$, its activity $a(B') := \sum_{b_i \in B'} a(b_i)$.
A territory $T$ is defined as a subset of the set of basic areas. The activity index of a territory $T$ is defined as 
$a(T) = \sum_{b_i \in T} a(b_i)$. The territory design problem is to assign each basic area to one of the ${k}$ territories, where $k > 1$ is the number of desired territories given in advance. A solution to the territory design problem is a set of territories $\mathcal{T} = \{T_1,\ldots , T_k\}$ such that for each pair of territories,  
$T_i \cap T_j = \emptyset$ and $T_1\cup\cdots\cup T_k=B$. A \emph{balance constraint} on the activity index demands that $a(T) \leq (1+\epsilon)\lceil\frac{a(B)}{k}\rceil \ \forall T \in \mathcal{T}$ for some \emph{imbalance tolerance} parameter $\epsilon \geq 0$. A solution is called \emph{feasible} if it obeys the balance constraint and other constraints placed on the \emph{contiguity} and \emph{compactness} of territories. 
The objective is to minimize the sum of pair-wise travel distances within the territories. Note, for small values of $\epsilon$ a feasible solution may not exist.

Consider an undirected graph $G=(V,E,c,\omega)$ 
with edge weights $\omega: E \to \MdR_{>0}$ and node weights
$c: V \to \MdR_{\geq 0}$.
Again, we extend $c$ and $\omega$ to sets, \ie , for a set $V' \subseteq V$ its node weight 
$c(V')\Is \sum_{v\in V'}c(v)$ and for a set $E' \subseteq E$ its edge weight $\omega(E')\Is \sum_{e\in E'}\omega(e)$.
The set $\Gamma(v)\Is \setGilt{u}{\set{v,u}\in E}$ denotes the neighbors of $v$.
A maximal connected component is a maximal subgraph of $G$ in which every pair of nodes are connected by a path. 
The \emph{graph partitioning problem} is looking for $k$ \emph{blocks} of nodes $V_1$,\ldots,$V_k$ 
that partition $V$, i.e., $V_1\cup\cdots\cup V_k=V$ and $V_i\cap V_j=\emptyset$
for $i\neq j$. Here, a \emph{balancing constraint} demands that 
$\forall i\in \{1,\ldots ,k\}\gilt c(V_i)\leq L_{\max}\Is (1+\epsilon)\lceil c(V)/k \rceil$ for
some imbalance parameter $\epsilon \geq 0$. 
Often the objective of graph partitioning problems is to minimize the total \emph{cut} $\sum_{i<j}w(E_{ij})$ where 
$E_{ij}\Is\setGilt{\set{u,v}\in E}{u\in V_i,v\in V_j}$. However, during the course of the paper we will modify this objective.
\csch{when the paper is finished, check whether we can remove definitions here}
\csch{define neighborhood graph}
\csch{define TSP tours etc...}

\subsection{Related Work}
There has been a \emph{huge} amount of research on territory design so that we refer the reader to the surveys \cite{duque2007supervised,kalcsics2005towards,zoltners20052004} for most of the material. For material related to graph partitioning, we refer the reader to the survey~\cite{SPPGPOverviewPaper}. 
Here, we focus on issues closely related to our main contributions. 

The basic formulation of territory design problems is due to Hess and Samuels~\cite{hess1971experiences} who used relaxations of integer linear programming models to tackle the problem.
Ronen \cite{ronen1983sales} presents a mixed integer formulation to minimize the total driving distance of salesmen. 
However,  it is assumed that each trip to a customer requires a separate trip which includes that the salesmen returns home, so that the average journey time is not measured. 

Using the same assumption, Zoltner and Sinha \cite{zoltners1983sales} develop a linear optimization model incorporating road networks for the distance computation to minimize the travel time. Additionally, the networks are used to represent connectivity between the territories which helped to compute contiguous territories with improved accessibility. Note that in both cases the objective also implicitly encodes compactness of the territories.

Forman and Yue \cite{Forman2003} use a genetic algorithm to compute congressional territories.
The basic idea is based on an encoding and on genetic operators that where originally used to solve the Traveling Salesman Problem.
Desired properties of the territories are integrated into the fitness function of the individuals. 
In this paper, we also use an evolutionary algorithm which instead encodes the individuals as partitions of a graph model.

Hess \etal \cite{hess1965} describe a mixed integer program based on the location-allocation approach. In each iteration basic areas are assigned to territory centers and afterwards the centers are updated. Compact territories are achieved by minimizing the sum of the distances of the basic areas to the territory centers. Since we also implement and extend this approach we go into more detail later.

Recently, Butsch \etal~\cite{bkn13} used a recursive geometric bipartitioning approach to assign basic areas to territories. The approach recursively computes bipartitions of the basic areas using their coordinates. However, this approach may compute longish, non-compact territories and has problems incorporating geographic obstacles such as rivers or mountains. We compare our algorithm against this approach in Section~\ref{s:experiments}.

\subsection{KaHIP}
\label{s:kaHIP}
Within this work, we use the open source multilevel graph partitioning framework KaHIP~\cite{kaffpa,kabapeE} (Karlsruhe High Quality Partitioning)\footnote{available at \url{http://algo2.iti.kit.edu/documents/kahip/} }.
More precisely, we modify the distributed evolutionary algorithm KaFFPaE contained therein to create partitions of our graph models which in turn yield solutions to the territory design problem. 
Hence, we shortly outline the main components of KaHIP. 

Besides the evolutionary algorithm, KaHIP implements many different algorithms, for example flow-based methods and more-localized local searches within a multilevel framework called KaFFPa, as well as several coarse-grained parallel and sequential meta-heuristics. 
The algorithms in KaHIP have been able to improve the best known partitioning results in the Walshaw Benchmark~\cite{soper2004combined} for many inputs using a short amount of time to create the partitions.
\paragraph{Evolutionary Algorithm Outline.} We now roughly outline the general structure of KaFFPaE since one of our algorithms employs a modified version of the evolutionary algorithm.
KaFFPaE starts with a population of individuals (in our case partitions of the graph) and evolves the population into different populations over several rounds. 
In each round, the evolutionary algorithm uses a selection rule based on the fitness of the individuals  of the population to select good individuals and combines them to obtain improved 
\begin{center}
\vspace*{-.8cm}
\begin{algorithm}[h]
\begin{algorithmic}
\normalsize
\STATE    create initial population $P$ 
\STATE    \textbf{while} stopping criterion not fulfilled 
\STATE    \quad \textit{select} parents $p_1, p_2$ from $P$ (using fitness function)
\STATE    \quad \textit{combine} $p_1$ with $p_2$ to create offspring $o$ (using a graph partitioner)
\STATE    \quad \textit{mutate} offspring $o$ 
\STATE    \quad \textit{evict} individual in population using $o$ 
\STATE    \textbf{return} the fittest individual that occurred
\end{algorithmic}
\caption{General Stucture of the Evolutionary Graph Partitioning Algorithm KaFFPaE.}
\label{alg:generalsteadystateEA}
\end{algorithm}
\end{center}
offspring. 
More precisely, KaFFPaE uses a tournament selection rule to select individuals for combination.
In KaFFPaE, the fitness function is set to the number of edges cut, however, we will modify the fitness function for the purpose of territory design. 
A combine operation then employs a graph partitioning framework to obtain an offspring having the good cuts of both input partitions. 
We refer the reader to \cite{kaffpaE} for more details.
Our algorithm generates one offspring per generation.  
The general structure of the evolutionary algorithm is depicted in Algorithm~\ref{alg:generalsteadystateEA}.

The algorithm is parallelized by giving each processing element its own population so that combine and mutation operations can be performed independently. 
This is combined with a scalable communication protocol to exchange high quality solutions between the processing elements over time. 
\section{Territory Design by Graph Partitioning}
\label{s:territorydesignbygp}

\label{s:graphpartitioning}
The territory design problem and the graph partitioning problem are closely related. 
Our approach to territory design using graph partitioning consists of two steps. 
At first we construct a graph that corresponds to the territory design problem. 
Then, a customized algorithm partitions the constructed graph into $k$ blocks. 
We develop a one-to-one mapping between basic areas and vertices in the constructed model.
Hence, a partition of our model yields a solution of the territory design problem.
Edges are defined by a \emph{neighboring} relation. 
We begin this section by explaining how we construct the graph to be partitioned, then define the fitness function that is used in the evolutionary algorithm and outline how everything is put together.

\subsection{Constructing the Graph}
\label{s:graphconstruction}
We now explain how we construct the graph that will be partitioned by the graph partitioning algorithm which in turn gives us a solution to the territory design problem.
Every basic area $b_i$ with activity $a_i$ corresponds to a node $v_i$ in the graph with weight $c_i := a_i$.
An edge between a node $v_i$ and $v_j$ exists if and only if the basic areas $b_i$ and $b_j$ are \emph{neighbors} as described below. 
In general, we set the weights of the edges in our model to one.  
However, as we will see later, the evolutionary graph partitioning algorithm that is used to partition the model takes the real distances into account to compute the value of the objective function of a solution. 

\paragraph{Edges in our Model.} The algorithm we use to compute the edges in our model is as follows. 
Note that we have multiple design goals for our model. 
First of all, two basic areas that are \emph{close} should be connected by an edge in the model.
Additionally, the model should not be too dense, \eg if the maximum degree of the graph is bounded, the graph should be connected and edges that are too long should also not be contained in the model.  
Roughly speaking, we perform two iterations of Kruskal's algorithm \cite{kruskal1956shortest} on a \emph{complete graph} where every node $v_i$ corresponds to a basic area $b_i$ and the edge weights are the distances between the basic areas $\omega(v_i,v_j) = d(b_i, b_j)$. Recall, that Kruskal's algorithm scans the edges of the graph in increasing order of their weight to grow a forest and that it adds the edges joining two trees in the forest.

After the first run of Kruskal's algorithm, \textit{every} edge that is in the minimum spanning tree (MST) computed by the algorithm is inserted into our model. 
For the second run of Kruskal's algorithm, we remove the just computed MST edges from the complete graph. 
Now, while the algorithm scans the edges in increasing order, it adds the current edge to our model if the maximum node degree in the current state of the model does not exceed a user defined parameter $\gamma > 1$ and if the length of the edge is smaller or equal to $\beta\cdot\omega_\text{avg}$. Here, $\beta$ is a user specified factor and $\omega_\text{avg}$ denotes the average edge weight of MST edges of the first iteration of the algorithm. 

\subsection{Fitness Function}
\label{s:objectivefunctiongraphpartitioning}

Recall that we are looking for a partition of the just defined graph model into $k$ blocks of nodes $P = \{V_{1} , \ldots , V_{k} \}$. 
Moreover, each subgraph induced by a block $V_{i}$ should to be \emph{contiguous}. 
In other words, if $n_{mcc}(V_i)$ is the number of maximal connected components in the subgraph induced by $V_i$, we want $n_{mcc}(V_i) = 1$ for each block $V_{i}$. 
We define the number of connected components of a partition $P$ as 
$n_{con}(P) = \sum_{i=1}^{k}{n_{mcc}(V_i)}$.
To guide the evolutionary algorithm towards contiguous blocks, we use a penalty approach.
More precisely, we use the factor $(1+ \alpha (n_{con}(P) - k))$ with penalty parameter $\alpha > 0$ in our objective function to ensure that non-contiguous territories are penalized. 
We then set the objective/fitness function of the evolutionary algorithm that partitions our model to 
\[ f_{obj}(P) = (1+ \alpha (n_{con}(P) - k))\sum_{i=1}^{k}{\sum_{u,v \in V_i }{{d(u,v)}}} \; , \]
where $d(u,v)$ is the average time needed to traverse the shortest route in the original network between the basic areas corresponding to $u$ and $v$ respectively.

\subsection{Overall Algorithm}
To tackle the territory design problem, we use the distributed evolutionary graph partitioner KaFFPaE to partition our model. 
We \emph{modify} the fitness function of the algorithm to the objective function presented in Section~\ref{s:objectivefunctiongraphpartitioning}.
Note that the partitioning algorithms within KaFFPaE still optimize the number of cut edges.
The amount of allowed imbalance is a constraint of the partition problem in KaFFPaE. 
Due to the local search algorithms in KaFFPaE no block of a partition will be empty.
However, it may be possible that partitions created during the course of the evolutionary algorithm are not contiguous.
Hence, whenever we create an individuum/partition, we try to make it contiguous.
This is done by grouping neighboring connected components of the partition so that each block becomes contiguous. 
More precisely, excess connected components are assigned to the neighboring block with the least activity.
Doing this may result in imbalanced partitions so that a rebalancing step is performed afterwards to ensure the balance constraint. We call the overall algorithm to tackle the territory design problem KaTeD (Karlsruhe Territory Design).
\vfill
\pagebreak

\section{Territory Design by Mixed Integer Programming}
\label{s:territorydesignbymip}

Besides the graph partitioning approach, we use a location-allocation method for territory design. In this section, we first describe the location-allocation method introduced by Hess \etal \cite{hess1965} and then present the details of our modified and extended location-allocation approach.

\subsection{The location-allocation approach by Hess \etal}
The idea of solving territory design problems by means of a location-allocation approach goes back to Hess \etal \cite{hess1965}. They apply the approach to the design of legislative districts. In this application, territories are supposed to be compact, contiguous, and balanced in terms of population.
Due to the complexity of the problem, Hess \etal \cite{hess1965} decompose it into two subproblems which are solved in an iterative manner:

\begin{itemize}
\item The \textbf{location} subproblem seeks to find (virtual) territory centers which are used to calculate the compactness measure. In the first iteration, the centers are a guess; in all subsequent iterations they are obtained by computing the center of gravity in each territory. 

\item In the \textbf{allocation} subproblem the basic areas are assigned to territory centers. To this purpose, Hess \etal \cite{hess1965} formulate an integer linear program.
Let $C$ denote the set of territory centers and $d_{b_i c_j}$ the distance between basic area $b_i \in B$ and center $c_j \in C$. Furthermore, define the following decision variables:
\begin{align*}
x_{b_i c_j} &= \begin{cases} 1 \text{\qquad if basic area $b_i$ is assigned to center $c_j$} \\ 0 \text{\qquad otherwise} \end{cases} \\
\end{align*}
Then, the model of Hess \etal \cite{hess1965} can be stated as follows:
\begin{alignat}{2}
&\text{min} \sum\limits_{b_i \in B} \sum\limits_{c_j \in C} d^2_{b_i c_j} a_i x_{b_i c_j} \label{obj_func} \\
&\sum\limits_{c_j \in C} x_{b_i c_j} = 1  \qquad & b_i \in B \label{constr_assign} \\
&(1-\epsilon) \frac{a(B)}{k} \leq \sum\limits_{b_i \in B} a_i x_{b_i c_j} \leq (1+\epsilon) \frac{a(B)}{k} \qquad & c_j \in C \label{constr_bal} \\
& x_{b_i c_j} \in \{0,1\} \qquad & b_i \in B, c_j \in C \label{constr_integr}
\end{alignat}
The objective function \eqref{obj_func} optimizes compactness which is measured as the sum of the squared distances between basic areas and associated territory centers, weighted by the basic areas' activity index. Constraints \eqref{constr_assign} in combination with the integrality conditions \eqref{constr_integr} ensure that each basic area is assigned to exactly one territory center. Balance is achieved by the constraints in \eqref{constr_bal} which guarantee that the activity index of a territory deviates by at most $\epsilon \cdot 100$ percent from the mean activity index.
Instead of solving the integer program, Hess \etal \cite{hess1965} set $\epsilon$ to zero, solve the linear programming relaxation and then resolve all fractional assignments.
\end{itemize}
Hess \etal \cite{hess1965} perform location and allocation alternately until the solution converges. The algorithm repeats if multiple initial guesses for the territory centers are available.

\subsection{Modifications and Extensions}

We now outline our modifications and extensions of the approach by Hess \etal \cite{hess1965}.
Roughly speaking, we use a well-known procedure to determine good initial centers, modify the balance constraint, solve the integer program directly and apply a multi-start procedure. We call our version of that approach KaLocAlloc (Karlsruhe Location Allocation).

\paragraph{Location Step.}
The quality of the solutions obtained by the location-allocation approach strongly depends on the selection of the initial centers, \ie , the centers used in the first iteration of the algorithm. A good initial set of centers should be well-distributed across the region under study. To this end, we adopt the initialization procedure from the k-means++ algorithm \cite{arthur2007}: Among all basic areas we pick the first center uniformly at random. Among all remaining basic areas the next center is picked with a probability which is proportional to the basic area's squared distance to the nearest center already chosen. This is repeated until ${k}$ centers have been selected.

\paragraph{Allocation Step.}
For the allocation of basic areas to centers we use model \eqref{obj_func} - \eqref{constr_integr} with one modification. In order to limit only the \textit{maximum} activity index of the territories, we drop the lower bound in constraints \eqref{constr_bal}:

\begin{alignat}{2}
&\sum\limits_{b_i \in B} a_i x_{b_i c_j} \leq (1+\epsilon) \frac{a(B)}{k} \qquad \forall \; c_j \in C \label{constr_bal2}
\end{alignat}

Although we drop the lower bound, solutions cannot contain empty territories. This is because the centers are picked among the basic areas and, therefore, at least the basic areas corresponding to centers are assigned to the respective center due to a distance of 0.
Another difference to Hess \etal \cite{hess1965} is that we do not solve the linear programming relaxation of the model, but the integer program.

\paragraph{Multi-start Procedure.}
As already mentioned, the selection of good initial centers is important to achieve solutions of high quality. Therefore, we apply a multi-start procedure in which the problem is solved multiple times with different initial centers at each start. The centers at each start are selected randomly according to the k-means++ scheme described above. After a user-defined number of starts the approach returns the best solution across all starts. Since in our case the objective is to minimize the sum of pairwise travel times between all basic areas of the same territory, the multi-start algorithm returns the best solution according to this criterion.
\vfill
\pagebreak
\section{Experimental Evaluation}
\label{s:experiments}

\paragraph{Methodology.}
All experiments were performed on real-world data provided by PTV Group. The test data comprises 15 instances whose sizes range from approximately 300 to 5,000 basic areas. The number of territories to be planned is given in the test data and varies from 3 to 46. Furthermore, road distances and travel times for the shortest path between each pair of basic areas were available and have been used for compactness evaluation.
We also compare our algorithms against the recursive partitioning approach by Butsch \etal\cite{bkn13} (BKNS) which has been provided by the authors. However, we perform only one repetition since the algorithm is deterministic.

We evaluated our approaches on all 15 instances with a time limit of 300 seconds. The imbalance tolerance parameter $\epsilon$ was set to 0.05 in all experiments. From the practical experience of PTV Group the chosen values for the time limit and for the imbalance tolerance are very acceptable values for human planners.

\paragraph{Parametrization.}
Parameters which are specific to territory design by graph partitioning were set as follows:
penalty parameter to increase connectedness $\alpha = 0.1$;
edge length factor in our model $\beta = 5$;
node degree bound in our model $\gamma = 20$.
We repeated both approaches 40 times using different random seeds for initialization and the average was taken.
Parameters which are specific to the location-allocation approach were set as follows: For the allocation step, the relative MIP optimality gap was set to 0.001. The maximum time spent on the allocation step was limited to 15 seconds, and the number of multi-starts was unrestricted. 

\paragraph{System.}
All experiments have been done on an Intel Xeon CPU E3-1245 at 3,4GHz having 16 GB RAM and Microsoft Windows 7.
The location-allocation approach was implemented in Java. We used Gurobi 5.6 to solve the integer linear program within the location-allocation approach.
The number of threads was set to 4 for KaTeD and KaLocAlloc to ensure comparability.

\subsection{Computational Results}
\begin{figure}
  \centering
   \includegraphics[width=0.75\textwidth]{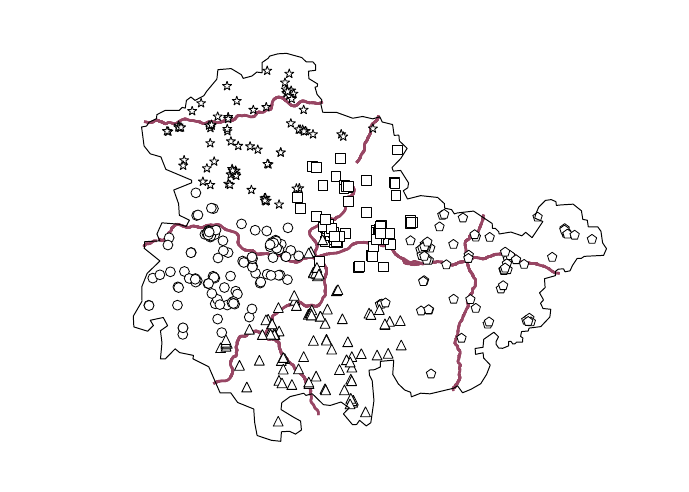}
  \caption{Visualization of the KaTeD result of instance BL16 (Thüringen, Germany).}
  \label{fig:thuringenkaffpaE}
\end{figure}

We shortly summarize the main results and present detailed per instance results in Table~\ref{table:results}.
First of all, the balance constraint is satisfied in all cases for all approaches.
In 8 of 15 instances, the average solution quality of KaTeD outperforms KaLocAlloc.
Considering the other 7 instances, the value of the objective of the territories computed by KaLocAlloc is 1.6 percent smaller than the ones computed by the evolutionary approach.
The recursive partitioning approach BKNS always yields worse results than both of our approaches, but is faster since we could perform only one repetition of the algorithm due to the fact that the algorithm is deterministic. On the largest instance, G02, the BKNS algorithm needed 203 seconds to compute a result. On average, BKNS yields 3.7 percent larger objectives than KaTeD and 4.3 percent larger objectives than KaLocAlloc.  The largest improvement over BKNS is obtained on instance BL16 (Thüringen) and amounts to 10.3 percent compared to the result of KaTeD which computes the best objective on that instance. 
Figures~\ref{fig:thuringenkaffpaE} to \ref{fig:BWMIP} compare the visual result computed by the different algorithms on two exemplary instances, Thüringen and Baden-Würtemberg in Germany.
\vfill
\section{Conclusion and Future Work}
\label{s:conclusion}
In this paper we addressed the territory design problem by developing graph theoretic models that also consider the underlying road network.  
The derived graph models enabled us to tackle the territory design problem by reducing it to a graph partitioning problem. The resulting graph partitioning problem is then solved by using a modified evolutionary graph partitioning algorithm takes the objective function of the territory design problem into account. 
On the other hand we extended an existing mixed integer programming formulation.
We tested and compared the algorithms on several real world instances. 

Important future work includes the integration of the objective function of the territory design problem directly into a multi-level graph partitioning algorithm. In particular, it would be interesting to define local search algorithms for our objective. On the other hand, it would be good to have a graph partitioning algorithm that can ensure connectedness of blocks.
%%=====================================================================
\paragraph*{Acknowledgements.} We would like to thank Alexander Butsch for providing us with an implementation of the recursive partitioning approach of Butsch \etal \cite{bkn13}.
%\section*{Bibliography} 
{
        \normalsize
\bibliographystyle{plain}
\bibliography{phdthesiscs,refs-parco,quellen_bender}
}
\vfill
\pagebreak

\begin{appendix}
\section{Tables and Pictures}
\begin{table}[H]
\caption{Sum of Pairwise Travel Times -- Average}
%\vspace*{.25cm}
\label{table:instances}
\centering
\normalsize
\begin{tabular}{ l | r | r | r | r | r }                   
Instance & \#Basic Areas & \#Territories & KaTeD avg & KaLocAlloc avg & BKNS avg\\
  \hline   
BL1  & \numprint{520}  & \numprint{6}  & \numprint{56080176}            &\textbf{\numprint{54157925}}   & \numprint{56943136}  \\ 
BL3  & \numprint{1265} & \numprint{13} & \textbf{\numprint{143300619}}  & \numprint{144435173}          & \numprint{146538257} \\ 
BL5  & \numprint{2440} & \numprint{20} & \numprint{239489621}           & \textbf{\numprint{238965327}} & \numprint{250332974} \\ 
BL6  & \numprint{514}  & \numprint{6}  &  \textbf{\numprint{109625833}}  & \numprint{110561195}     & \numprint{118085932} \\ 
BL7  & \numprint{1038} & \numprint{9}  &   \textbf{\numprint{48151371}}   & \numprint{48456559}          & \numprint{50664649}  \\ 
BL8  & \numprint{1252} & \numprint{11} & \textbf{\numprint{159771590}}  & \numprint{160239521}          & \numprint{172659265} \\ 
BL9  & \numprint{2019} & \numprint{18} & \numprint{260762258}           & \textbf{\numprint{258027628}} & \numprint{271895453} \\ 
BL11 & \numprint{427}  & \numprint{3}  & \numprint{28990949}            & \textbf{\numprint{28796550}}  & \numprint{29745445}  \\ 
BL13 & \numprint{289}  & \numprint{3}  & \numprint{50774924}            & \textbf{\numprint{50326699}}  & \numprint{50301370}  \\ 
BL14 & \numprint{493}  & \numprint{5}  & \textbf{\numprint{55345856}}   & \numprint{55723702}           & \numprint{55844700}  \\ 
BL15 & \numprint{284}  & \numprint{3}  & \numprint{41953610}            & \textbf{\numprint{40997882}}  & \numprint{41426146}  \\ 
BL16 & \numprint{428}  & \numprint{5}  & \textbf{\numprint{44396659}}   & \numprint{44564384}           & \numprint{48977328}  \\ 
G01  & \numprint{4472} & \numprint{41} & \numprint{469430937}           & \textbf{\numprint{459630225}} & \numprint{479442432} \\ 
G02  & \numprint{4971} & \numprint{45} & \textbf{\numprint{597515028}}  & \numprint{591514229}          & \numprint{632623485} \\ 
G03  & \numprint{2241} & \numprint{20} & \textbf{\numprint{333612833}}  & \numprint{329314024}          & \numprint{343170613} \\ 
\end{tabular}
\end{table}

\begin{table}[H]
\caption{Sum of Pairwise Travel Times -- Minimum and Maximum}
\label{table:results}
\centering
\normalsize
\begin{tabular}{ l | r | r | r | r | r}                   
Instance & KaTeD min & KaLocAlloc min & KaTeD max & KaLocAlloc max & BKNS min/max\\
  \hline   
BL1  & \numprint{54918613}           & \textbf{\numprint{53663743}}  & \numprint{56276617}            &\textbf{\numprint{54819237}}   & \numprint{56943136} \\ 
BL3  & \textbf{\numprint{142464954}} & \numprint{143406994}          & \textbf{\numprint{143688813}}  & \numprint{145504176}          & \numprint{146538257}\\ 
BL5  & \numprint{235752526}          & \textbf{\numprint{233110996}} & \textbf{\numprint{244426840}}  & \numprint{245775302}          & \numprint{250332974}\\ 
BL6  &    \textbf{\numprint{109140977}} & \numprint{109548283}   & \textbf{\numprint{109893677}}  & \numprint{111336978}          & \numprint{118085932}\\ 
BL7  & \textbf{\numprint{47118442}}  & \numprint{48293411}               & \numprint{48591642}            & \textbf{\numprint{48569098}}  & \numprint{50664649} \\ 
BL8  & \numprint{158613738}          & \textbf{\numprint{158257253}} & \textbf{\numprint{160899255}}  & \numprint{161325525}          & \numprint{172659265}\\ 
BL9  & \numprint{257234587}          & \textbf{\numprint{253678835}} & \numprint{262697498}           & \textbf{\numprint{261954880}} & \numprint{271895453}\\ 
BL11 & \numprint{28973547}           & \textbf{\numprint{28796550}}  & \numprint{28992672}            & \textbf{\numprint{28796550}}  & \numprint{29745445} \\ 
BL13 & \numprint{50774924}           & \textbf{\numprint{50217885}}  & \numprint{50774924}            & \textbf{\numprint{50446097}}  & \numprint{50301370} \\ 
BL14 & \textbf{\numprint{54102563}}  & \numprint{55456757}           & \numprint{57347598}            & \textbf{\numprint{55810459}}  & \numprint{55844700} \\ 
BL15 & \numprint{41255274}           & \textbf{\numprint{40987865}}  & \numprint{41971516}            & \textbf{\numprint{41011058}}  & \numprint{41426146} \\ 
BL16 & \textbf{\numprint{44083530}}  & \numprint{44455849}           & \numprint{44609070}            & \textbf{\numprint{44631114}}  & \numprint{48977328} \\ 
G01  & \numprint{459200690}          & \textbf{\numprint{447119646}} & \numprint{473914701}           & \textbf{\numprint{473030108}} & \numprint{479442432}\\ 
G02  & \numprint{581718616}          & \textbf{\numprint{575445420}} & \numprint{612582152}           & \textbf{\numprint{606800497}} & \numprint{632623485}\\ 
G03  & \numprint{329746198}          & \textbf{\numprint{326456733}} & \numprint{336518355}           & \textbf{\numprint{333657431}} & \numprint{343170613}\\ 
\end{tabular}
\end{table}

\begin{figure}
  \centering
   \includegraphics[width=0.75\textwidth]{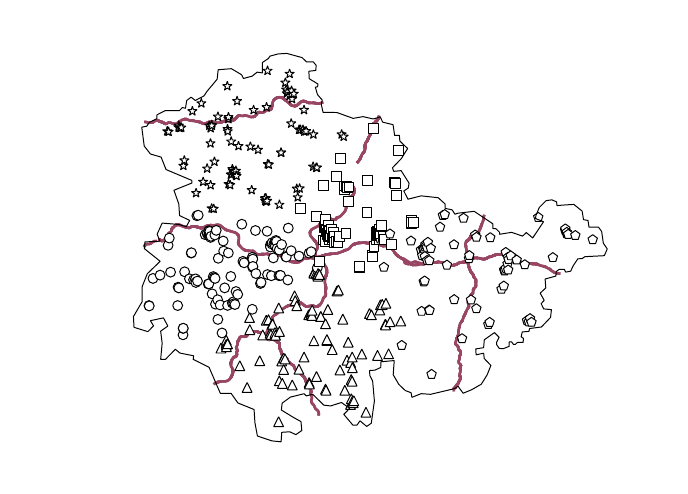}
  \caption{Visualization of the  KaLocAlloc result of instance BL16 (Thüringen, Germany).}
  \label{fig:thuringenMIP}
\end{figure}

\begin{figure}
  \centering
   \includegraphics[width=0.75\textwidth]{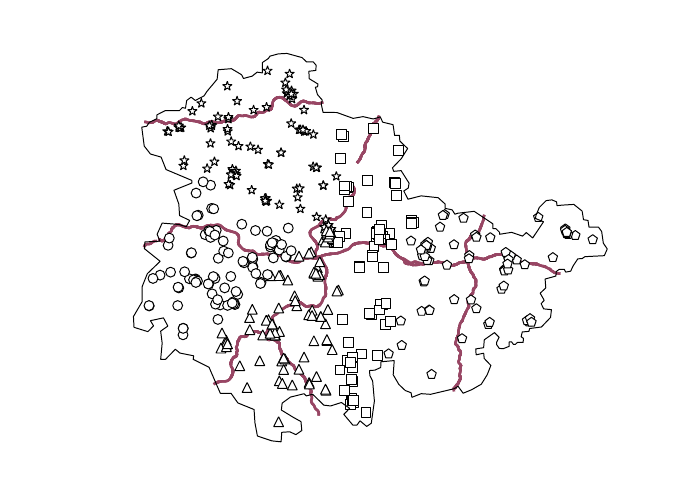}
  \caption{Visualization of the  BKNS result of instance BL16 (Thüringen, Germany).}
  \label{fig:thuringenMIP}
\end{figure}

\begin{figure}
  \centering
   \includegraphics[width=.5\textwidth]{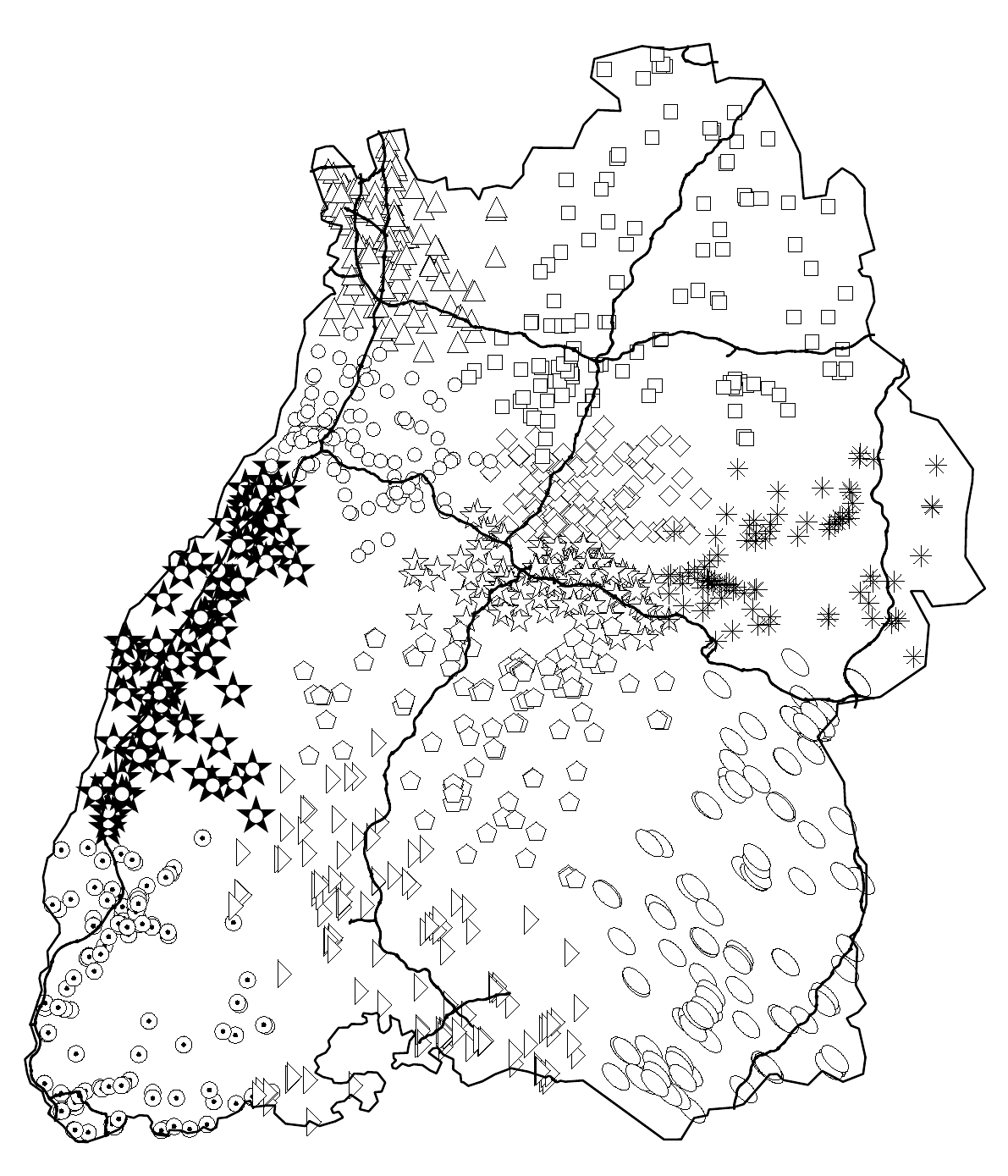}
  \caption{Visualization of the KaTeD result of instance BL8 (Baden-Wuerttemberg, Germany)}
  \label{fig:BWKaFFPaE}
\end{figure}

\begin{figure}
  \centering
  \includegraphics[width=.5\textwidth]{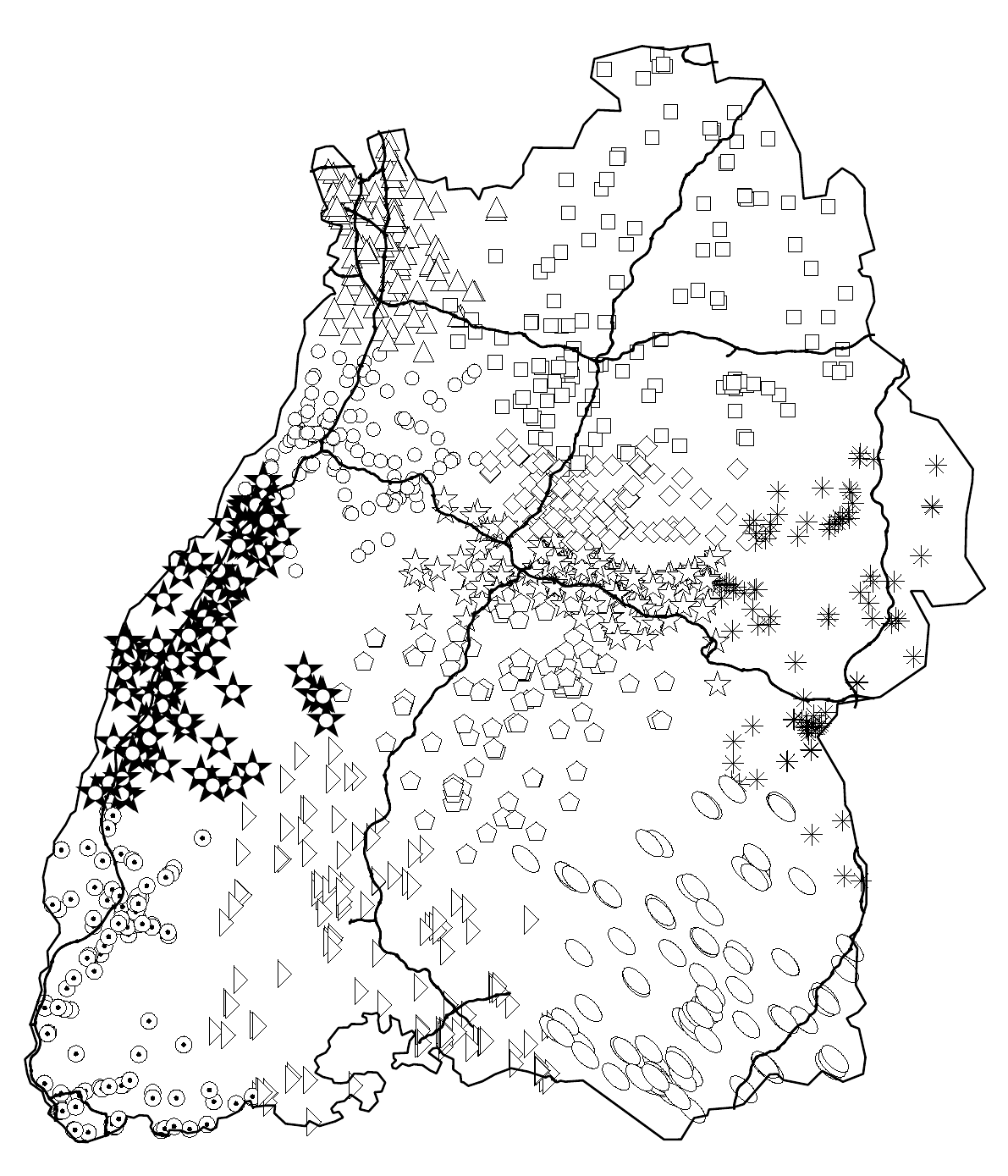}
  \caption{Visualization of the KaLocAlloc result of instance BL8 (Baden-Wuerttemberg, Germany).}
  \label{fig:BWMIP}
\end{figure}

\begin{figure}
  \centering
  \includegraphics[width=.5\textwidth]{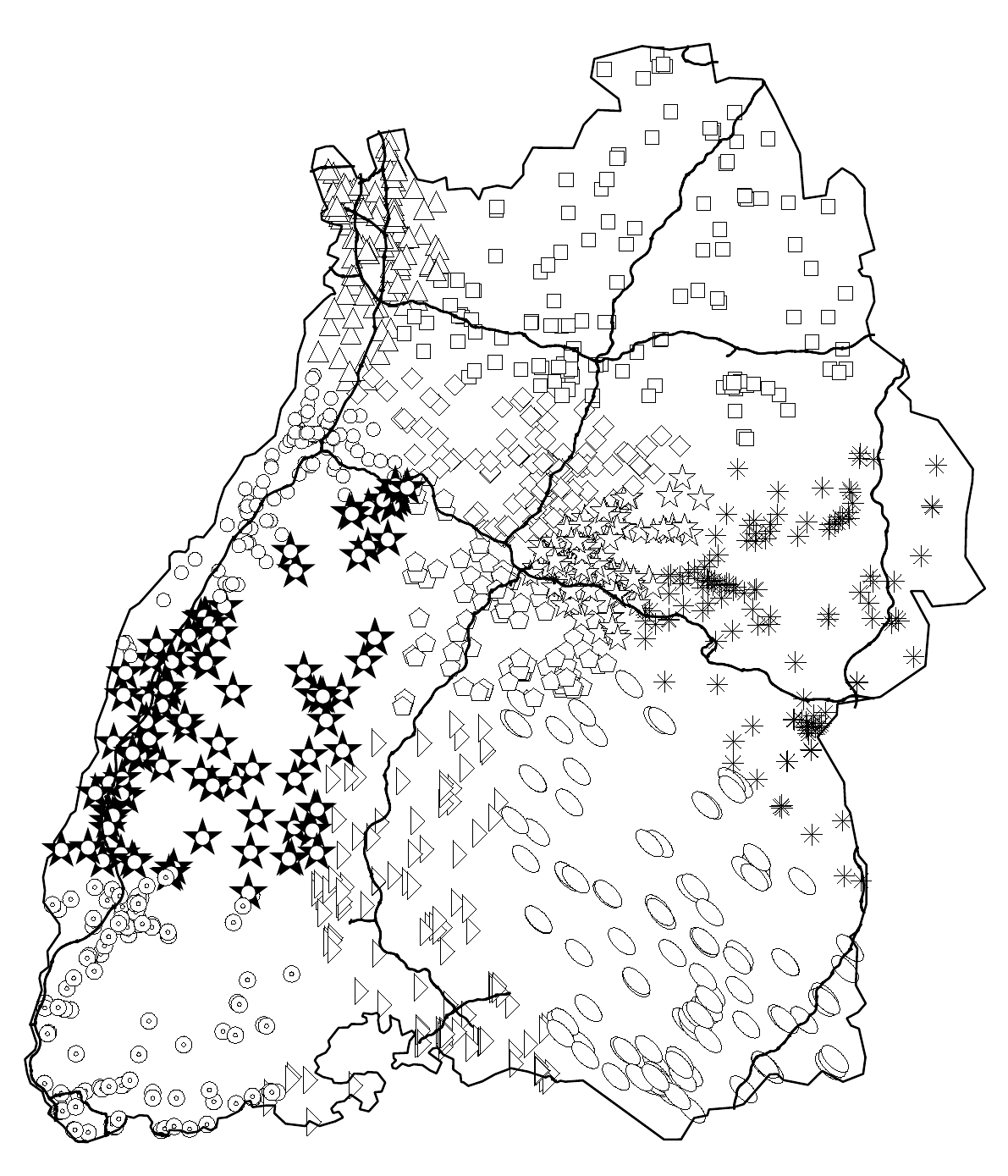}
  \caption{Visualization of the BKNS result of instance BL8 (Baden-Wuerttemberg, Germany).}
  \label{fig:BWMIP}
\end{figure}

\end{appendix}
\vfill
\end{document}